\documentclass[12pt]{article}

\usepackage{latexsym}
\usepackage{calc}
\usepackage{datetime}
\usepackage{amssymb, amsmath, amsfonts, listings, mathrsfs} 
\usepackage{xcolor} 
\usepackage[normalem]{ulem} 

\usepackage{comment} 
\usepackage{float}

\usepackage{algorithmic}
\usepackage[algoruled,boxed,lined]{algorithm2e}

\usepackage{graphicx}
\usepackage[labelfont=bf]{caption}

\newcounter{hours}\newcounter{minutes}

\newcommand{\printtime}{
\setcounter{hours}{\time/60}
\setcounter{minutes}{\time - \value{hours}*60}
\thehours h \theminutes m}

\newcommand{\stkout}[1]{\ifmmode\text{\sout{\ensuremath{#1}}}\else\sout{#1}\fi} 

\def\nr{\par \noindent}

\def\Def{\stackrel{\mathrm{def}}{=}}

\def\inter{{\rm int \,}}

\def\dom{{\rm dom \,}}
\def\beq{\begin{equation}}
\def\eeq{\end{equation}}

\def\R{\mathbb{R}}
\def\E{\mathbb{E}}

\def\BI{\begin{itemize}}
\def\EI{\end{itemize}}

\newcommand{\SetEQ}{\setcounter{equation}{0}}
\newcommand{\refLE}[1]{\ensuremath{\stackrel{(\ref{#1})}{\leq}}}
\newcommand{\refEQ}[1]{\ensuremath{\stackrel{(\ref{#1})}{=}}}

\newcommand{\refLEI}[2]{\ensuremath{\stackrel{(\ref{#1})_{#2}}{\leq}}}
\newcommand{\refEQI}[2]{\ensuremath{\stackrel{(\ref{#1})_{#2}}{=}}}

\newtheorem{theorem}{Theorem}
\newtheorem{lemma}{Lemma}
\newtheorem{corollary}{Corollary}

\newtheorem{assumption}{Assumption}
\newtheorem{definition}{Definition}

\newtheorem{example}{Example}
\newtheorem{remark}{Remark}
\newcommand{\proof}{\bf Proof: \rm \nr}
\newcommand{\qed}{\hfill $\Box$ \nr \medskip}
\newcommand{\half}{\mbox{${1 \over 2}$}}

\def\ba{\begin{array}}
\def\ea{\end{array}}
\def\beann{\begin{eqnarray*}}
\def\eeann{\end{eqnarray*}}
\def\bea{\begin{eqnarray}}
\def\eea{\end{eqnarray}}

\def\BT{\begin{theorem}}
\def\ET{\end{theorem}}
\def\BL{\begin{lemma}}
\def\EL{\end{lemma}}
\def\BC{\begin{corollary}}
\def\EC{\end{corollary}}
\def\BE{\begin{example}}
\def\EE{\end{example}}
\def\BD{\begin{definition}}
\def\ED{\end{definition}}
\def\BR{\begin{remark}}
\def\ER{\end{remark}}
\def\BAS{\begin{assumption}}
\def\EAS{\end{assumption}}
\def\BI{\begin{itemize}}
\def\EI{\end{itemize}}

\def\BMP{\begin{minipage}{9.5cm}}
\def\EMP{\end{minipage}}
\def\MPT{\begin{minipage}{11.5cm}}
\def\EPT{\end{minipage}}

\def\la{\langle}
\def\ra{\rangle}

\def\QF{\hspace{5ex} \Box}
\def\QR{\hfill \Box}

\title{
Theorem of Alternative for Extended Homogeneous Linear System and its Application in Conic Optimization}
\author{Yurii Nesterov\thanks{
Emeritus professor at CORE/INMA, UClouvain, Belgium. Email: Yurii.Nesterov@uclouvain.be. 
\newline Emeritus professor at Corvinus University of Budapest. Fractional professor at School of Data Sciences\\ (Chinese University of Hong Kong, Shenzhen)
} 
}

\date{\vspace{3ex} \normalsize 
March 22, 2026\\ 
\vspace{1ex} [$\;$version 4.0, file: {FeasIPM$\backslash$InfeasIPM4.tex}$\;$]
\\ \vspace{2ex} Printout: \printtime, \today
}

\begin{document}
\maketitle

\abstract{In this paper, we develop a new framework for constructing infeasible-start primal-dual methods for Conic Optimization. Our approach can be seen as a straightforward consequence of Gordan Theorem of Alternative. Given by the target upper bound $\epsilon > 0$ for the duality gap as the only input parameter, we form an auxiliary convex  problem of minimizing barrier function with linear equality constraints. Its solution can be easily transformed to the requested output. This function can be minimized by different schemes of Unconstrained Optimization, with possible quadratic convergence in the end of the process. In our paper, we analyze the Damped Newton Method and a short-step path-following scheme. For both of them, we prove polynomial-time complexity results. Our methods are able to benefit from the hot-start opportunities. We can ensure the residual of the  linear equality constraints in the primal and dual problems to be at the level of machine accuracy, independently on the accuracy parameter $\epsilon$.}

\vspace{10ex}\noindent
{\bf Keywords:} Conic Optimization, Primal-Dual Problem, Interior-Point Methods, Self-Concordant Barriers, Infeasible-Start Methods.

\thispagestyle{empty}

\newpage\setcounter{page}{1}

\section{Introduction}
\setcounter{equation}{0}

\vspace{1ex}\noindent
{\bf Motivation.} As it is widely recognized now, theorems of alternatives play an important role in Optimization Theory (e.g. \cite{Broyden}). One of the first theorems of this type was proved in 1873 by Gordan \cite{Gordan} for homogeneous systems of linear inequalities coupled with a positive orthant. Later on in 1902, Farkas \cite{Farkas} published a non-homogeneous version of this result, which became one of the fundamentals of Convex Analysis (see \cite{Rock}).

These theorems are usually applied to resolve some questions related to existence and boundedness of  optimal sets of different optimization problems. Hence, in the more recent years, there are still publications proposing new variants of these results (e.g. \cite{Olvi}) and finding new ways for their justifications (e.g. \cite{Mil}). Further developments are also focused on extensions of the initial theorems onto general nonlinear cones (e.g. \cite{SKAR}).

In this paper, we present a variant of Gordan Theorem for general cones, where a part of variables is free. As we will see, this is exactly the format suitable for development of {\em Infeasible-Start} primal-dual Interior-Point Methods (IPMs) for general convex cones. Such a development is the main goal of this paper.

Infeasible start IPMs are the most powerful modern methods for solving the problems of Conic Optimization. They can start from an arbitrary interior point inside the convex cone and generate a sequence of iterates, for which both the duality gap and the residual in the primal-dual system of linear equality constraints are vanishing. For that, they introduce additional projective variables (one or two) and let the minimization sequence go to infinity. If the problem is feasible, then the projective variable goes to infinity too. The final output is obtained by dividing the current primal-dual point by one of the projective variables. This operation explains the main drawback of this approach: in the output,  both the accuracy in function value and the residual of linear equality constraints are always of the same order. However, in many practical problems we do not need an extremely high accuracy in the function value. At the same time, very often, the linear constraints must be satisfied with the machine accuracy.

Another drawback of these two approaches is impossibility to use the advantages of hot start. Even if our primal-dual starting point is close somehow to the solution, the global complexity bounds of these schemes remain unchanged. In our paper, we are going to address both these issues.

In order to position properly our results, let us explain two existing homogenization approaches for constructing infeasible start IPMs (see \cite{Tuncel} for reviewing other techniques).
Both of them are applied to the following primal-dual pair of Conic Problems:
\beq\label{prob-PD0}
\ba{rcl}
f^* & = & \min\limits_{x \in K} \{ \la c, x \ra: \; A x = b \} \; = \; \max\limits_{y \in \E_y, \, s \in K^*} \{ \la b, y \ra: \; s + A^* y = c \},
\ea
\eeq
where $c \in \E^*_x$, $b \in \E^*_y$, and $A: \E_x \to \E_y^*$ (see Section \ref{sc-Not} for notation). The optimal solution $(x^*,s^*,y^*)$ is characterized by the following conditions:
\beq\label{eq-Opt0}
\ba{rclrcl}
 s^* + A^* y^* & = & c, & A x^* & = & b,\\
 \\
 \la c, x^* \ra - \la b, y^* \ra & = & 0, & x & \in & K, \quad s \in K^*.
\ea
\eeq

In the first approach \cite{NestINF}, we introduce a single projective variable $t \in \R$ and form the following {\em shifted} feasible set:
\beq\label{eq-OptS1}
{\cal F}_1 = \left\{(x,s,y,t):
\ba{rclrcl}
 s + A^* y & = & \bar s + t c, & A x& = & A \bar x + tb,\\

 \la c, x \ra - \la b, y\ra & = & \la c, \bar x \ra, & x & \in & K, \quad s \in K^*
\ea \right\}.
\eeq
where $\bar x \in \inter K$ and $\bar s \in \inter K^*$ are the initial reference points in the primal-dual cone. It is clear that if $(\hat x, \hat s, \hat y, \hat t)$ is a recession direction of the set ${\cal F}$ with $\hat t > 0$, then the point $(\hat x /\hat t, \hat s/ \hat t, \hat y / \hat t)$ satisfies conditions (\ref{eq-Opt0}). On the other hand, it is well known that the recession directions of an unbounded convex set can be found by minimizing its self-concordant barrier. This simple observation leads to several infeasible-start IPMs \cite{NestINF}. All of them  solve the following optimization problem:
\beq\label{prob-M1}
\inf\limits_{(x,s,y,t) \in {\cal F}_1}  \Big\{ F(x) + F_*(s) \Big\}.
\eeq

The second approach is based on so-called {\em self-dual embedding} (see \cite{YTM} for the initial version for Linear Optimization and \cite{AY,DKRT,LSZ,PS} for further extensions). In this framework, in the case of general convex cones, it is necessary to find a recession direction of the following set \cite{NYT}:
\beq\label{eq-OptS2}
\ba{c}
{\cal F}_2 = \Big\{\{(x,s,y,t,\kappa):\;  s + A^* y = \bar s + A^* \bar y - \bar t c + t c, \quad A x = A \bar x - \bar t b + tb,\\
 \la c, x \ra - \la b, y\ra + \kappa = \bar \kappa+ \la c, \bar x \ra - \la b, \bar y \ra,\\
x \in  K, \quad s \in K^*, \quad t \geq 0, \quad \kappa \geq 0  \Big\},
\ea
\eeq
where $\bar x \in \inter K$, $\bar s \in \inter K^*$, $\bar y \in \R^m$, $\bar t > 0$, and $\bar \kappa > 0$ are the reference points in the corresponding sets.
For finding this direction, we can solve the following minimization problem:
\beq\label{prob-M2}
\inf\limits_{(x,s,y,t,\kappa) \in {\cal F}_2}  \Big\{ F(x) + F_*(s) - \ln t - \ln \kappa\Big\},
\eeq
which is below unbounded provided that our initial problem (\ref{prob-PD0}) has nonempty interior. Same as for problem (\ref{prob-M1}), the output is generated by dividing the last iterate by the projective variable $t > 0$. Hence, the accuracy in linear equality constraints is proportional to that of duality gap.

Our approach can be seen as a mixture of these two approaches. We {\em derive} it from one fundamental fact in Convex Analysis, a version of Gordan Theorem of Alternative for extended linear systems and convex cones.

\vspace{1ex}\noindent
{\bf Contents.} The paper is organized as follows. In Section \ref{sc-Not}, we introduce notation and present some general results used in the paper. In Section \ref{sc-Prob}, we prove a modified version of Gordan Theorem of Alternative for an extended system of linear equality constraints and a general convex cone. This theorem is used in Section \ref{sc-Opt} in order to derive a convex optimization problem (\ref{prob-Shift}) with equality constraints. Its solution can be transformed to a strictly feasible $\epsilon$-solution of the initial problem. In all our methods, accuracy $\epsilon > 0$ plays a role of input parameter.

In Section \ref{sc-Newt}, we analyze behavior of Damped Newton Method as applied to problem~(\ref{prob-Shift}). For this method, we prove a polynomial-time complexity bound of the order $O(\nu \ln {1 \over \epsilon})$, where $\nu$ is the parameter of self-concordant barrier for the feasible cone.

In Section \ref{sc-Short}, we analyze performance of a short-step path-following scheme. We show that it achieves the region of quadratic convergence in $O(\nu^{1/2} \ln {1 \over \epsilon})$ iterations. In Section~\ref{sc-Hot}, we show that our methods benefit from using a hot-start reference points.

We finish the paper with Conclusion in Section \ref{sc-Conc} and Appendix in Section \ref{sc-App}, where we put some discussion and auxiliary results on self-concordant functions.

\section{Notations and Generalities}\label{sc-Not}
\SetEQ

Let $\E$ be a finite-dimensional space and $\E^*$ be the space of linear functions on $\E$. For $x \in \E$ and $s \in \E^*$, denote by $\la s, x \ra$ the value of function $s$ at $x$. We use the same notation $\la \cdot, \cdot \ra$ for different spaces. Thus, its actual meaning is defined by the context.
For a linear operator $A: \E_1 \to \E^*_2$, we denote by $A^*$ its {\em adjoint operator}:
\beq\label{def-AD}
\ba{rcl}
\la A x, y \ra & = & \la A^* y, x \ra, \quad x \in \E_1, \; y \in  \E_2.
\ea
\eeq
Thus, $A^*: \E_2 \to \E_1^*$. The operator $A: \E \to \E^*$ is positive semidefinite (notation $A \succeq 0$) if $\la A x, x \ra \geq 0$ for all $x \in \E$. We write $A \succeq B$ if $A-B \succeq 0$. 

In case of $\E = \R^n$, we have $\E^* = \R^n$ and
$$
\ba{rcl}
\la s, x \ra & = & \sum\limits_{i=1}^n s^{(i)} x^{(i)}, \quad \| x \| \; \Def \; \la x, x \ra^{1/2}, \quad s, x \in \R^n.
\ea
$$
Vector $e \in \R^n$ is the vector of all ones, and $e_i$, $i = 1, \dots , n$, are coordinate vectors in $\R^n$. 

Sometimes it is convenient to treat any vector $c \in \E^*$ as a linear operator acting from $\R$ to $\E^*$ is accordance to the following rule:
$$
\ba{rcl}
c \tau & = & \tau \cdot c, \quad \tau \in \R.
\ea
$$
Then $c^*: \E \to \R$ and $c^* x = \la c, x \ra$ for $x \in \E$.

For function $f(\cdot)$ with open domain $\dom f \subseteq \E$, we denote by
$$
\ba{rcl}
D^kf(x)[h_1, \dots, h_k], \quad x \in \dom f, \; h_i \in \E, \; i = 1, \dots, k,
\ea
$$
its $k$th directional derivative at $x$ along $k$ directions. It is a symmetric $k$-linear form. If all directions are the same, we use notation $D^kf(x)[h]^k$. Thus, under the mild assumptions, for the gradient and the Hessian of function $f(\cdot)$, we have the following relations:
$$
\ba{rcl}
Df(x)[h] & = & \la \nabla f(x), h \ra, \quad D^2f(x)[h]^2 \; = \; \la \nabla^2 f(x) h, h \ra.
\ea
$$

Let us recall some facts from the theory of self-concordant functions (e.g. \cite{LN}). 

\BD 
Let closed convex function $f \in {\cal C}^3(\E)$ have open domain. It is called {\em self-concordant} (SCF) if its third directional derivative is bounded by an appropriate power of the second one: 
\beq\label{def-SCF}
\ba{rcl}
D^3f(x)[h]^3 & \leq & 2 \Big( D^2 f(x)[h]^2 \Big)^{3/2}, \quad x \in \dom f, \; h \in \E.
\ea
\eeq
\ED

If $\dom f$ contains no straight line, then its Hessian is positive definite at any $x \in \dom f$. In this paper, we mainly deal with such functions. However, $\dom f$ still can be unbounded. In this case, we have the following Theorem on Recession Direction:
\beq\label{eq-RecDir}
\ba{rcl}
\la \nabla^2 f(x) d, d \ra^{1/2} & \leq & - \la \nabla f(x), d \ra 
\ea
\eeq
for any $x \in \inter \dom f$ and any {\em recession direction} $d$ of the set $\dom f$.

For any SCF, its dual function $f_*(s) = \max\limits_{x \in \dom f} [ - \la s, x \ra - f(x)]$ is also self-concordant. Note that $- \nabla f(x) \in \dom f_*$ for all $x \in \dom f$, and
$$
\ba{rcl}
\nabla f_*(s) & = & - x(s) \; \Def\; - \arg\max\limits_{x \in \dom f} [ - \la s, x \ra - f(x) ], \quad s \in \dom f_*
\ea
$$
For any $x \in \dom f$ and $s \in \dom f_*$, we have
\beq\label{eq-FDual}
\ba{rclrcl}
\nabla f(-\nabla f_*(s)) & = & -s, & \nabla^2 f(-\nabla f_*(s)) & = & [ \nabla^2 f_*(s)]^{-1},\\
\\
\nabla f_*(-\nabla f(x)) & = & -x &
\nabla^2 f_*(-\nabla f(x)) & = & [ \nabla^2 f(x)]^{-1}.
\ea
\eeq

One of the main properties of SCF is that the {\em Dikin Ellipsoid} defined as
$$
\ba{rcl}
W^f_r(x) & = & \{ u \in \E:\; \la \nabla^2 f(x)(y-x),y - x \ra \leq r^2\}, \quad x \in \dom f,
\ea
$$
belongs to $\dom f$ for any $r \in [0,1)$. 

In this paper, we often use the local norms defined by the Hessians of SCF. For $x \in \dom f$, $s \in \dom f_*$, $u \in \E$, and $g \in \E^*$, we adopt the following notation:
$$
\ba{rclrcl}
\| u \|_{\nabla^2 f(x)} & = & \la \nabla^2 f(x) u, u \ra^{1/2}, &  \quad \| g \|_{\nabla^2 f(x)} & = & \la  g, [\nabla^2 f(x)]^{-1} g \ra^{1/2}, \\
\\
\| g \|_{\nabla^2 f_*(s)} & = & \la g, \nabla^2 f_*(s) g \ra^{1/2}, &  \quad \| u \|_{\nabla^2 f_*(s)} & = & \la [\nabla^2 f_*(s)]^{-1} u, u \ra^{1/2}.
\ea
$$
If no ambiguity arise, we use the corresponding shortcuts $\| u \|_x$, $\| g \|_x$, $\| g \|_s$, and $\| u \|_s$. Hence, the sense of notation $\| a \|_b$ depends on the sets and the spaces containing $a$ and $b$. 

Let us recall now the properties of {\em self-concordant barriers } (SCB) for convex cones.
\BD\label{def-SCB}
A self-concordant function $F(\cdot)$ is called a {\em self-concordant barrier}, if there exists a constant 
\beq\label{eq-DefNu}
\ba{rcl}
\nu & \geq & 1
\ea
\eeq
such that for all $x \in \dom F$ and $h \in \E$ we have
\beq\label{eq-SCB}
\ba{rcl}
\la \nabla F(x), h \ra^2 & \leq & \nu \la \nabla^2 F(x) h, h \ra.
\ea
\eeq
The constant $\nu$ is called the {\em parameter} of the barrier.
If $\nabla^2 F(x) \succ 0$, then the condition (\ref{eq-SCB}) is equivalent to the following:
\beq\label{def-SCB1}
\ba{rcl}
\lambda_F^2(x) \; \Def \; \| \nabla F(x) \|_x^2 & \leq & \nu, \quad x \in \dom F.
\ea
\eeq
\ED
This value is responsible for complexity of the set $\dom F$ for corresponding IPMs. 

In the theory of IPMs, the most important feasible sets are convex cones. A cone $K$ is called {\em proper} if it is closed convex and pointed (contains no straight lines). For such cones, SCBs possess a natural property of {\em logarithmic homogeneity}:
\beq\label{def-HomB}
\ba{rcl}
F(\tau x) & \equiv & F(x) - \nu \ln \tau, \quad x \in \inter K, \; \tau > 0.
\ea
\eeq
This identity has several important consequences: for all $x \in \inter K$ and $\tau > 0$, we have
\beq\label{eq-Hom12}
\ba{rclrcl}
\nabla F(\tau x) & = & {1 \over \tau} \nabla F(x), & D^k F(\tau x) & = & {1 \over \tau^k} D^k F(x), \; k \geq 2,
\ea
\eeq
\beq\label{eq-HomFX}
\ba{rcl}
\la \nabla F(x), x \ra & = & - \nu, \quad \la \nabla^2 F(x) x, x \ra \; = \; \nu, \quad \| \nabla F(x) \|^2_x \; = \; \nu.
\ea
\eeq
\beq\label{eq-Hom1}
\ba{rcl}
\nabla^2 F(x) x & = & - \nabla F(x).
\ea
\eeq
Thus, the parameter of a logarithmically homogeneous SCB for a proper cone (which we call {\em regular barrier}) is equal to the degree of logarithmic homogeneity.

It is important that the dual barrier $F_*(\cdot)$ for the regular barrier $F(\cdot)$ is a regular barrier for the dual cone
$$
\ba{rcl}
K^* & = & \Big\{ s \in \E^*: \; \la s, x \ra \geq 0, \; \forall x \in K \Big\},
\ea
$$
with the same value of barrier parameter. 

\section{Theorem of Alternative for Convex Cones}\label{sc-Prob}
\SetEQ

Let $K_z \subset \E_z$ be a closed convex pointed cone with nonempty interior. For two spaces $\E_u$ and $\E_y$, consider linear operators $B: \E_y \to \E_u^*$, and  $C: \E_z \to \E_u^*$. Then, we have
$$
\ba{rcl}
B^*: \; \E_u \to \E_y^*, \quad C^*: \E_u \to \E_z^*.
\ea
$$
We assume that operator $(C,B): \E_z \times \E_y \to \E_u^*$ has full row rank.

Consider the following pair of convex sets:
$$
\ba{rcl}
{\cal F}_p & = & \{ (z,y): \; z \in \inter K_z, \; y \in \E_y: \; C z + B y = 0\},\\
\\
{\cal F}_d & = & \{ u \in \E_u:\; C^* u \in K^*_z, \; B^* u = 0 \}.
\ea
$$

Let us prove the following variant of Gordan Theorem.
\BT\label{th-Alt}
Cone ${\cal F}_p$ is not empty if and only if ${\cal F}_d\equiv \{0\}$.
\ET
\proof
Let $(z, y) \in {\cal F}_p$ and $u \in {\cal F}_d$, $u \neq 0$. Then,
$$
\ba{rcl}
0 & < & \la C^* u , z\ra \; = \; \la C z , u \ra \; = \; - \la B y, u \ra \; = \; - \la B^* u, y \ra \; = \; 0,
\ea
$$
and this is impossible.

Assume now that ${\cal F}_d \equiv \{0\}$. Let us choose an arbitrary point $\bar s \in \inter K^*_z$ and consider the set $Q = \{ u \in \E_u: \; \bar s + C^* u \in K^*_z, \; B^* u = 0 \}$. In view of our assumption, this set has no recession directions. Hence, it is bounded and its relative interior is non-empty. 

Let us consider an arbitrary self-concordant barrier $F(\cdot)$ for the cone $K_z$.
Then, there exists a unique solution to the following optimization problem:
\beq\label{prob-Prim}
\ba{rcl}
\min\limits_{u \in Q} F_*(\bar s + C^* u).
\ea
\eeq
Denote by $y_*$ the optimal Lagrange multipliers for the equality constraints of this problem, and by $u_*$ its optimal solution. Then, the first-order optimality conditions tell us that
$$
\ba{rcl}
C \nabla F_*(\bar s + C^* u_*) & = & B y_*.
\ea
$$
Since $z_* \Def - \nabla F_*(\bar s + C^* u_*) \in \inter K^*_z$, we get $(z_*,y_*) \in {\cal F}_p$.
\qed
\BC\label{cor-Alt}
Either ${\cal F}_p \neq \emptyset$ or ${\cal F}_d \neq \{0\}$. Situation ${\cal F}_p \neq \emptyset$ and ${\cal F}_d \neq \{0\}$ is impossible.
\EC
\proof
Indeed, either ${\cal F}_d = \{0\}$, or ${\cal F}_d \neq \{0\}$. These events are clearly complementary.
\qed

\BC\label{cor-Alt2}
We never have at the same time ${\cal F}_p =\emptyset$ and ${\cal F}_d = \{0\}$. 
\EC

It is important that Theorem \ref{th-Alt} provides us with a constructive way of finding an element from the set ${\cal F}_p$.  For that, we need to solve the optimization problem (\ref{prob-Prim}). In the subsequent sections, we will show how this result can be used for constructing Infeasible-Start IPMs for solving a primal-dual pair of Conic Optimization Problems.

\section{Infeasible-Start Formulation for \\ Conic Optimization}\label{sc-Opt}
\SetEQ

Let $K_x \subset \E_x$ be a proper cone. Consider the following primal-dial pair of Conic Problems:
\beq\label{prob-PD}
\ba{rcl}
f^* & = & \min\limits_{x \in K} \{ \la c, x \ra: \; A x = b \} \; \geq \; f_*\; = \; \max\limits_{y \in \E_y, \, s \in K^*} \{ \la b, y \ra: \; s + A^* y = c \},
\ea
\eeq
where $c \in \E^*_x$, $b \in \E^*_y$, and $A: \E_x \to \E_y^*$. We assume that $A$ has full row rank and $c \not\in \rm{Im}(A^*)$. Under the following strict feasibility assumption:
\beq\label{ass-Feas}
\ba{rcl}
\hat x & \in & \inter K_x, \; A \hat x \; = \; b, \quad \hat s \; \in \inter K^*_x, \quad \hat s + A^* \hat y \; = \; c,
\ea
\eeq
both problems are solvable and $f^* = f_*$ (see, for example, \cite{LongStep}).

Let $K_x$ admit a $\nu$-self-concordant logarithmically homogeneous barrier $F(\cdot)$. Then we can define the primal and dual central paths as follows:
$$
\ba{rcl}
x(t) & = & \arg \min\limits_x \{ t \la c, x \ra +F(x): \; A x = b \}, \quad t > 0,\\
\\
(s(t),y(t)) & = & \arg\min\limits_{s,y} \{ - t \la b, y \ra + F_*(s): \; s + A^* y = c \}.
\ea
$$
It can be proved that the components of the primal and dual central paths satisfy the following relations (e.g. \cite{LongStep}):
\beq\label{eq-RCP}
\ba{rcl}
s(t) & = & - {1 \over t} \nabla F(x(t)),\\
\\
\la s(t), x(t) \ra & = & \la c, x(t) \ra - \la b, y(t) \ra \; = \; {\nu \over t},\\
\\
F(x(t)) + F_*(s(t) )& = & - \nu + \nu \ln t.
\ea
\eeq

The problem (\ref{prob-PD}) can be solved by IPMs based on a path-following strategy. However, for that it is necessary to compute first a strictly feasible primal-dual point. It appears that Theorem \ref{th-Alt} gives us a possibility to avoid this potentially expensive computation.

Note that the optimal primal dual solution $(x^*,s^*,y^*)$ of problem (\ref{prob-PD}) is characterized by the following conditions:
\beq\label{eq-Opt}
\ba{c}
A x^* = b, \quad s^* + A^* y = c,\quad \la c,x^*\ra - \la b,y^* \ra = 0,\\
\\
x^* \in K_x, \quad s^* \in K^*_x.
\ea
\eeq
However, it is impossible to find a solution to this system by interior-point methods since the equality constraints there do not intersect the interior of the feasible cones.

Let us make the formulation of our problem suitable for applying Theorem \ref{th-Alt}. For that, let us add to the system (\ref{eq-Opt}) a projective variable $t \geq  0$ and introduce an accuracy parameter $\epsilon > 0$. This gives us the following set:
\beq\label{def-FP}
\ba{rcl}
{\cal F}_p & =& \{ z = (x,s,t) \in K_z, \; y \in \E_y: \; A x = t b, \; s + A^* y = t c, \; \la c,x \ra - \la b, y \ra = t \epsilon \}\\
\\
K_z & = & \inter K_x \times \inter K_x^* \times \inter \R_{+}.
\ea
\eeq
Thus, in the notation of Theorem \ref{th-Alt}, we have $E_z = \E_x \times \E^*_x \times \R$, $\E^*_u = \E_x^* \times \E_y^* \times  \R$, and
$$
\ba{rcl}
C & = & \left( \ba{rrr} 0 & I & -c \\  -A & 0 & b \\ c^* & 0 & - \epsilon \ea \right): \; \E_z \to \E_u^* , \quad B \; = \; \left( \ba{c} A^* \\  0 \\ - b^* \ea \right): \E_y \to \E^*_u.
\ea
$$
In this setting, it is clear that the operator $(C,B)$ has full row rank.

In view of assumption (\ref{ass-Feas}), ${\cal F}_p \neq \emptyset$. Therefore, by Theorem \ref{th-Alt}, we have
$$
\ba{rcl}
{\cal F}_d & \Def & \{ u \in \E_u: \; C^*u \in K^*_z, \; B^* u = 0 \} \; \equiv \; \{ 0 \}.
\ea
$$

Coming back to our initial notation and decomposing the variable $u$ as $u = (x,y,\tau)$, we get the following representation of the set ${\cal F}_d$:
\beq\label{def-FD}
\ba{r}
{\cal F}_d = \Big\{ (x,y,\tau) \in \E_u: \;  \tau c - A^*y \in K^*_x, \; x \in K_x, \\
- \la c, x \ra + \la b, y \ra - \epsilon \tau \geq 0, \; A x = \tau b \Big\}.
\ea
\eeq
Since $0 \leq \la  \tau c - A^*y,x \ra = \tau [ \la c, x \ra - \la b,y \ra] = \tau [ \la c, x \ra - \la b,y \ra + \epsilon] - \epsilon \tau^2$, we could add to ${\cal F}_d$ a redundant constraint $\tau \leq 0$.

Thus, our alternatives look now as follows (see Corollary \ref{cor-Alt}):
\beq\label{eq-AltCP}
\ba{l}
\mbox{\em Either there exist points $x_* \in \inter K_x$ and $y_* \in \E_y$ such that}\\
c- A^* y_* \in \inter K^*, \quad A x_* = b, \quad \la c, x_* \ra - \la b, y_* \ra = \epsilon,\\
\mbox{\em or,}\\
\mbox{\em There exist points $\hat x \in K_x$ and $\hat y \in \E_y$ with $\hat \tau \leq 0$ such that}\\
\hat \tau c - A^* \hat y \in K_x^*, \quad A \hat x = \hat \tau b, \quad \la c, \hat x \ra - \la b, \hat y\ra + \hat \tau \epsilon \leq 0,
\ea
\eeq
and these situations are complementary. By homogeneity, the second alternative can be reduced to two cases: $\hat \tau = 0$ and $\hat \tau = -1$,  which are not mutually exclusive. In the later case, the linear inequality is redundant.

Let us fix some reference points $\bar x \in \inter K$, $\bar s \in \inter K^*$, and $\bar \tau > 0$. We assume that ${\cal F}_p \neq \emptyset$. Then ${\cal F}_d$ contains only the origin, and the following problem is solvable:
\beq\label{prob-Shift}
\ba{c}
\min\limits_{u = (x,y,\tau): \, \atop A x = \tau b} \Big\{\Phi(u) \Def F(\bar x + x) + F_*( \bar s + \tau c - A^*y) - \ln( \bar \tau  - \la c, x \ra + \la b, y \ra - \epsilon \tau)  \Big\},
\ea
\eeq
where $\Phi(\cdot)$ is a self-concordant barrier for its domain with parameter $\bar \nu = 2 \nu + 1$.

From Theorem \ref{th-Alt}, we know that the optimal solution of this problem provides us with a feasible solution of the set ${\cal F}_p$ defined by (\ref{def-FP}). Let us confirm this theoretical prediction by a direct computation.

Denote by $u_* = (x_*,y_*, \tau_*)$ the optimal solution of problem (\ref{prob-Shift}) and by $\Phi_*$ its optimal value. Let $\lambda \in \E_y$ be the optimal  Lagrange multipliers for equality constraints in (\ref{prob-Shift}). Define the following Lagrange function:
$$
\ba{rcl}
{\cal L}(x,y,\tau,\lambda) & = & F(\bar x + x) + F_*( \bar s + \tau c - A^*y) - \ln( \bar \tau  - \la c, x \ra + \la b, y \ra - \epsilon \tau) \\
\\
& & + \la \tau b - A x, \lambda \ra.
\ea
$$
The first-order stationarity conditions for this function look as follows:  $A x_* = \tau_* b$ and
\beq\label{eq-FOrder0}
\ba{rcl}
\nabla F(\bar x+x_*) + {1 \over \omega_*} c - A^* \lambda_* & = & 0,\\
\\
- A \nabla F_*(\bar s + \tau_* c - A^*y_*) - {1 \over \omega_*} b & = & 0,\\
\\
\la c, \nabla F_*(\bar s + \tau_* c - A^*y_*) \ra + {1 \over \omega_*} \epsilon + \la b, \lambda_* \ra & = & 0,
\ea
\eeq
where $\lambda_*$ is the vector of optimal Lagrange multipliers for the equality constraints, and 
$$
\ba{rcl}
\omega_* & = & \bar \tau  - \la c, x_*\ra + \la b, y_* \ra - \epsilon \tau_* \; > \; 0.
\ea
$$
Hence, defining
\beq\label{def-PDSol}
\ba{rcl}
x_{\epsilon} & = & - \omega_* \nabla F_*(\bar s + \tau_* c - A^*y_*), \quad s_{\epsilon} \; = \; - \omega_*  \nabla F(\bar x+x_*), \quad y_{\epsilon} \; = \; \omega_* \lambda_*,
\ea
\eeq
we get a strictly feasible primal-dual pair of problem (\ref{prob-PD}) with the duality gap
$$
\ba{rcl}
\la s_{\epsilon}, x_{\epsilon} \ra & = & \la c, x_{\epsilon} \ra - \la b, y_{\epsilon} \ra \; \refEQI{eq-FOrder0}{3} \;  \epsilon.
\ea
$$
Thus, the point $u_*$ allows us to {\em compute} the feasible point from ${\cal F}_p$.

\section{Infeasible Start Newton Method}\label{sc-Newt}
\SetEQ

The simplest way of solving the problem 
\beq\label{prob-Phi}
\ba{rcl}
\Phi_* & = & \min\limits_{u = (x,y,\tau)} \Big\{ \Phi(u): \; A x = \tau b \Big\}
\ea
\eeq
is to apply the Damped Newton Method. It is important that we already have a strictly feasible starting point $0 \in \dom \Phi$. Let $u_*$ be the optimal solution of problem (\ref{prob-Phi}).

Denote by $\lambda(u)$ the restricted local norm of the gradient of function $\Phi(\cdot)$ at $u$:
\beq\label{def-GNorm}
\ba{rcl}
\lambda^2(u) & \Def & \max\limits_{h = (h_x,h_y,h_{\tau}) \in \E_u} \Big\{ 2 \la \nabla \Phi(u), h \ra - \la \nabla^2 \Phi(u)h,h \ra: \; A h_x = h_{\tau} b \Big\}.
\ea
\eeq
And let $h(u)$ be the optimal solution of this problem. Then the scheme of the Damped Newton Method is as follows:
\beq\label{eq-DNM}
\mbox{\fbox{$\;u_0 = 0, \quad u_{k+1} = u_k - {h(u_k) \over 1 + \lambda(u_k)}, \quad k \geq 0\;$}}
\eeq

It is well known that outside the region of quadratic convergence, this method decreases the value of the objective function by an absolute constant. Thus, it needs at most 
$$
O\Big(\Phi(0)- \Phi_*\Big)
$$
iterations for entering the region of quadratic convergence. Hence, we need to find some upper bounds for the initial gap in function value.

Let us start from the choice of reference points. We choose
\beq\label{eq-Choice}
\ba{rcl}
\bar \tau & = & 1, \quad \bar s \; = \; - \nabla F(\bar x),
\ea
\eeq
where $\bar x$ is an arbitrary reference point in $\inter K_x$. Then
\beq\label{eq-Phi0}
\ba{rcl}
\Phi(0) & = & F(\bar x) + F_*(-\nabla F(\bar x)) \; \refEQI{eq-HomFX}{1} \; - \nu.
\ea
\eeq

It remains to estimate from below the value $\Phi_*$. Let us represent it in the dual form.
$$
\ba{rcl}
\Phi_* & = & \min\limits_u \max\limits_{\tilde x, \tilde s,t,\lambda } \Big\{ - \la t \tilde s, \bar x + x \ra - F_*(t\tilde s) - \la \bar s + \tau c - A^* y, t \tilde x \ra - F(t \tilde x) \\
\\
& &- (\bar \tau - \la c, x \ra + \la b, y \ra - \epsilon \tau)t + \ln t + 1  + \la \tau b - A x, t \lambda\ra \Big\}\\
\\
& \refEQ{def-HomB} &  \max\limits_{\tilde x, \tilde s,t, \lambda} \min\limits_u \Big\{ - t [\bar \tau + \la \tilde s, \bar x \ra + \la \bar s, \tilde x \ra ] - F_*(\tilde s)  - F(\tilde x) + t  \la c - \tilde s - A^* \lambda , x \ra \\
\\
& & + t \la A \tilde x - b, y \ra + t\tau [ \epsilon - \la c, \tilde x \ra + \la b, \lambda \ra] + \bar \nu \ln t + 1\Big\}.
\ea
$$
Thus,
\beq\label{eq-RepPhi}
\ba{rclr}
\Phi_* & = & 1+ \max\limits_{\tilde x, \tilde s,t, \lambda} \Big\{  \bar \nu\ln t  - t [\bar \tau + \la \tilde s, \bar x \ra + \la \bar s, \tilde x \ra ] - F_*(\tilde s)  - F(\tilde x):\\
\\
& &  \hspace{15ex} c = \tilde s + A^* \lambda,\;  A \tilde x = b, \;   \la c, \tilde x \ra - \la b, \lambda \ra = \epsilon \Big\}\\
\\
& = & - 2\nu+ \max\limits_{(\tilde x, \tilde s, \lambda) \in {\cal F}} \Big\{   - \bar \nu \ln {\bar \tau + \la \tilde s, \bar x \ra + \la \bar s, \tilde x \ra \over \bar \nu} - F_*(\tilde s)  - F(\tilde x): \\
& & \hspace{18ex} \la c, \tilde x \ra - \la b, \lambda \ra = \epsilon  \Big\},
\ea
\eeq
where ${\cal F} \Def \Big\{ (x,s,y):\;  A x = b,\; s + A^* y = c, \; x \in K_x, \; s \in K^*_x \Big\}$.

In order to get a lower bound for $\Phi_*$, we need to estimate the size of the primal-dual central path with respect to the reference point. For that, let us define the minimal $\sigma \geq 1$ such that
\beq\label{def-Sigma}
\ba{c}
{1 \over \sigma} \bar x \; \preceq_{K_x} \; x(1) \; \preceq_{K_x} \; \sigma \bar x, \quad {1 \over \sigma} \bar s\; \preceq_{K^*_x} \;  s(1) \; \preceq_{K^*_x} \; \sigma \bar s.
\ea
\eeq
Then, for any $(\tilde x, \tilde s) \in K_x \times K^*_x$, we have $\la \tilde s, \bar x \ra + \la \bar s, \tilde x \ra \leq \sigma[ \la \tilde s, x(1) \ra + \la s(1), \tilde x \ra]$. If we choose $(\tilde x, \tilde s, \lambda) \in {\cal F}$, then
$$
\ba{rcl}
\la s(1) - \tilde s, x(1) - \tilde x \ra = 0.
\ea
$$
If in addition, we have $\la \tilde s, \tilde x \ra = \epsilon$, then we obtain the following relation:
\beq\label{eq-NE}
\ba{rcl}
 \la \tilde s, x(1) \ra + \la s(1), \tilde x \ra & = &  \la s(1), x(1) \ra + \la \tilde s, \tilde x \ra  \; \refEQI{eq-RCP}{1} \; \nu + \epsilon.
 \ea
 \eeq
 
Let us choose now $t_{\epsilon} = {\nu \over \epsilon}$. Then $\la s(t_{\epsilon}), x(t_{\epsilon}) \ra = \la c, x(t_{\epsilon}) \ra - \la b, y(t_{\epsilon}) \ra \refEQI{eq-RCP}{2} \epsilon$, and in accordance to representation (\ref{eq-RepPhi})$_2$, we have
$$
\ba{rcl}
\Phi(0) - \Phi_* & \refLE{eq-Phi0} &  \nu +  \bar \nu \ln {1 + \la s(t_{\epsilon}), \bar x \ra + \la \bar s, x(t_{\epsilon}) \ra \over \bar \nu } + F_*(s(t_{\epsilon}))  + F(x(t_{\epsilon}))\\
\\
& \leq & \nu +  \bar \nu \ln {1 + \sigma [ \la s(t_{\epsilon}), x(1) \ra + \la s(1), x(t_{\epsilon}) \ra] \over \bar \nu} + F_*(s(t_{\epsilon}))  + F(x(t_{\epsilon}))\\
\\
& \refEQ{eq-NE}  &  \nu +  \bar \nu \ln {1 + \sigma [ \nu + \epsilon] \over \bar \nu} + F_*(s(t_{\epsilon}))  + F(x(t_{\epsilon}))\\
\\
& \refEQI{eq-RCP}{3} &  \bar \nu \ln {1 + \sigma [ \nu + \epsilon] \over \bar \nu}  + \nu \ln {\nu \over \epsilon}.
\ea
$$

Thus, we have proved the following theorem.
\BT\label{th-Newt}
Let $\sigma$ be defined by (\ref{def-Sigma}) and the reference points are chosen by (\ref{eq-Choice}). Then
method (\ref{eq-DNM}) enters the region of quadratic convergence at most in
\beq\label{eq-UpDNM}
\ba{c}
O \Big( \bar \nu \ln {1 + \sigma [ \nu + \epsilon] \over \bar \nu }  + \nu \ln {\nu \over \epsilon}\Big)
\ea
\eeq
iterations.
\ET

Note that the estimate (\ref{eq-UpDNM}) does not depend on the quality of the starting point. Instead, it depends logarithmically on a much less sensitive characteristic, the choice of the reference points.

\section{Infeasible-Start Short-Step Path-Following \\ Scheme}\label{sc-Short}
\SetEQ

Let us show how we can solve the problem (\ref{prob-Phi}) by a short-step path-following scheme. 

For that, we need to explain first the procedure of choosing the reference points. We start from the point $\tilde w = (\tilde x, \tilde s, \tilde y)$ such that
$$
\ba{rcl}
\tilde x & \in \inter K, \quad \tilde s \; \in \; \inter K^*.
\ea
$$
These points reflect somehow our information about $\epsilon$-solution of problem (\ref{prob-PD}). 

After that, we define 
\beq\label{eq-Rule}
\ba{rclrcl}
\bar x & = & - \nabla F_*(\tilde s), & \bar s & = & - \nabla F(\tilde x), \quad \bar \tau \; = \; 1.
\ea
\eeq
If we choose, for example $\tilde y = 0$ and $\tilde s = - \nabla F(\tilde x)$, then we come to the rule (\ref{eq-Choice}). Now we can define the parameter $\sigma \geq 1$ by the conditions (\ref{def-Sigma}).

Define now the vector $g = (g_x,g_y,g_{\tau}) \in \E_x^* \times \E_y^* \times \R$ with the following components:
\beq\label{eq-ComP}
\ba{rcl}
g_x & = & \tilde s + A^* \tilde y - c, \quad g_y \; = \; b - A \tilde x , \quad g_{\tau} \; = \; \la c, \tilde x \ra - \la b, \tilde y \ra - \epsilon.
\ea
\eeq
This gives us a possibility to introduce the  following {\em dual central path}, which corresponds to the target vector $g$:
\beq\label{def-HotCP}
\ba{rcl}
u_d(t) & \Def & \arg\min\limits_{u = (x,y,\tau)} \Big\{ \Phi_t(u) \Def t \la g, u \ra + \Phi(u): \; A x = \tau b \; \Big\}, \; 0 \leq t \leq 1.
\ea
\eeq
The rationality behind all these definitions is clear from the following lemma.
\BL\label{lm-Dual}
For the dual path defined by (\ref{eq-Rule})-(\ref{def-HotCP}), we have
$u_d(0) = u_*$ and $u_d(1) = 0$.
\EL
\proof
Since $\Phi_0(\cdot) \equiv \Phi(\cdot)$, the first equality is trivial. Further, we have
$$
\ba{rcl}
{\partial \Phi_1(0) \over \partial x} & = & \nabla F(\bar x) + c + g_x \; = \; A^* \tilde  y,\\
\\
{\partial \Phi_1(0) \over \partial y} & = & - A \nabla F_*(\bar s)  - b + g_y \; = \; 0,\\
\\
{\partial \Phi_1(0) \over \partial \tau} & = &  \la c, \nabla F_*(\bar s) \ra + \epsilon + g_{\tau} \; = \; - \la b, \tilde y \ra.
\ea
$$
This means that $u_d(1) = 0$.
\qed

Thus, we can use the starting point $u=0$ in order to follow the dual path $u_d(t)$ 
as $t \to 0$ (see Section~5.3.4 in \cite{LN}). 
Let us define the restricted local norm of vector $v \in \E_u^*$:
\beq\label{def-GNormR}
\ba{rcl}
\lambda^2_u(v) & \Def & \max\limits_{h = (h_x,h_y,h_{\tau}) \in \E_u} \Big\{ 2 \la v, h \ra - \la \nabla^2 \Phi(u)h,h \ra: \; A h_x = h_{\tau} b \Big\},
\ea
\eeq
and let $h_u(v)$ be the optimal solution of this problem.
Then, we choose $u_0 = 0$, $t_0 = 1$, and iterate
\beq\label{met-SSPF}
\fbox{$\ba{rclrcl}
t_{k+1} & = & \left( t_k - {\hat \gamma \over \lambda_{u_k}(g)} \right)_+, & u_{k+1} & = & u_k - {1 \over 1 + \xi_k} h_{u_k}(\nabla \Phi_{t_{k+1}}(u_k))
\ea$}
\eeq
where $
\xi_k = {\lambda_k^2 \over 1 + \lambda_k}$ and $\lambda_k = \lambda_{u_k} (\nabla \Phi_{t_{k+1}}(u_k))$. It can be proved that with parameters $\hat \beta = 0.126$ and $\hat \gamma = 0.164$, we have 
\beq\label{eq-RateSS}
\ba{c}
\lambda_{u_k}(\nabla \Phi_{t_{k}}(u_k)) \; \leq \; \hat \beta,\quad
t_{k+1} \; \leq \; \left( 1 - {\hat \gamma \over \hat \beta + \sqrt{\bar \nu}} \right) t_k, \quad k \geq 0.
\ea
\eeq

Thus, method (\ref{met-SSPF}) has global linear rate of convergence. However, in a small neighborhood of the solution $u_*$, we can ensure the local quadratic convergence by applying a standard Newton Method. Let us derive the corresponding switching condition and estimate the number of iterations of method (\ref{met-SSPF}) required to meet it.

By Theorem 5.2.2 in \cite{LN}, the Damped Newton Method 
\beq\label{met-SNewton}
\ba{rcl}
u_{k+1}& = & u_k - {h_{u_k}(\nabla \Phi(u_k)) \over 1 + \lambda_{u_k}(\nabla \Phi(u_k))} 
\ea
\eeq
has the following local rate of convergence:
\beq\label{eq-RateDNM}
\ba{rcl}
\lambda_{u_{k+1}}(\nabla \Phi(u_{k+1})) & \leq & 2 \lambda^2_{u_k}(\nabla \Phi(u_k)).
\ea
\eeq
Thus, the region of quadratic convergence of this method is described by inequality:
\beq\label{eq-RQC}
\ba{rcl}
\lambda_{u}(\nabla \Phi(u)) & < & \half.
\ea
\eeq
Since
$$
\ba{rcl}
\lambda_{u_k}(\nabla \Phi(u_k)) & \leq & \lambda_{u_k}(\nabla \Phi_{t_k}(u_k)) + t_k \lambda_{u_k}(g),
\ea
$$
from (\ref{eq-RateSS})$_1$, we conclude that the switching rule for method (\ref{met-SSPF}) can be as follows:
\beq\label{eq-Switch}
\ba{rcl}
t_k \lambda_{u_k}(g) & < & \tilde \beta \Def \half - \hat \beta.
\ea
\eeq
Thus, we need to find an upper bound for the value 
$\lambda_{u_k}(g)$.
For that, we need to justify some lower bounds for the matrices $\nabla^2 \Phi(u_k)$.

Let us establish first some properties of the optimization problem (\ref{prob-Shift}) written in an extended form:
\beq\label{prob-ShiftX}
\ba{c}
\min\limits_{\hat u = (x,s,y,\tau)} \Big\{ \hat \Phi(\hat u) \Def F(\bar x + x) + F_*(s) - \ln (\bar \tau - \la c, x \ra + \la b, y \ra - \epsilon \tau): \hat u \in \hat Q \Big\},
\ea
\eeq
where 
$$
\ba{rr}
\hat Q \Def \Big\{ \hat u = (x,s,y,\tau):  & \bar x + x \in \inter K, \; s \in \inter K^*, \; 
s = \bar s + \tau c -A^* y, \quad \\
& A x = \tau b, \quad \bar \tau > \la c, x \ra - \la b, y \ra + \epsilon \tau \;\Big\}.
\ea
$$

First of all, let us show that the set $\hat Q$ is bounded.
\BL\label{lm-QBound}
For any $\hat u = (x,s,y,\tau) \in \hat Q$, we have $\bar x + x \in K$, $s \in K^*$, and
\beq\label{eq-QBound}
\ba{rcl}
\la s(t_{\epsilon}), \bar x + x \ra + \la s, x(t_{\epsilon}) \ra & < & \bar \tau + \la \bar s, x(t_{\epsilon}) \ra + \la s(t_{\epsilon}), \bar x \ra,
\ea
\eeq
where $t_{\epsilon} = {\nu \over \epsilon}$. Moreover,
\beq\label{eq-GapUp}
\ba{rcl}
\bar \tau - \la c, x \ra + \la b, y \ra - \epsilon \tau & \leq  & 
\bar \tau + \la \bar s, x(t_{\epsilon}) \ra + \la s(t_{\epsilon}), \bar x \ra.
\ea
\eeq
\EL
\proof
Let $\hat u = (x,s,y,\tau) \in \hat Q$. Then
\beq\label{eq-Gap}
\ba{rcl}
\la c, x \ra - \la b, y \ra & = & \la c - A^* y(t_{\epsilon}), x \ra + \la A x, y(t_{\epsilon}) \ra - \la b, y \ra\\
\\
& = & \la s(t_{\epsilon}), x \ra + \tau \la b, y(t_{\epsilon}) \ra - \la A^*y, x(t_{\epsilon}) \ra\\
\\
& = & \la s(t_{\epsilon}), x \ra + \tau \la b, y(t_{\epsilon}) \ra - \la\bar s + \tau c - s, x(t_{\epsilon}) \ra\\
\\
& \refEQI{eq-RCP}{2} & \la s(t_{\epsilon}), x \ra + \la s, x(t_{\epsilon}) \ra - \tau {\nu \over t_{\epsilon}} - \la \bar s, x(t_{\epsilon}) \ra\\
\\
& = & \la s(t_{\epsilon}), \bar x + x \ra + \la s, x(t_{\epsilon}) \ra - \tau \epsilon - \la s(t_{\epsilon}), \bar x \ra - \la \bar s, x(t_{\epsilon}) \ra.
\ea
\eeq
Thus, 
$$
\ba{rcl}
\bar \tau & > & \la c, x \ra - \la b, y \ra + \epsilon \tau \; \refEQ{eq-Gap} \;
\la s(t_{\epsilon}), \bar x + x \ra + \la s, x(t_{\epsilon}) \ra - \la s(t_{\epsilon}), \bar x \ra - \la \bar s, x(t_{\epsilon}) \ra,
\ea
$$
and this is (\ref{eq-QBound}). Finally, in view of representation (\ref{eq-Gap}), we have
$$
\ba{rcl}
\bar \tau - \la c, x \ra + \la b, y \ra - \epsilon \tau & \leq  & 
\bar \tau + \la s(t_{\epsilon}), \bar x \ra + \la \bar s, x(t_{\epsilon}) \ra.\QF
\ea
$$

Let use now the following consequences of equations (\ref{eq-SLow}), (\ref{eq-XLow}):
\beq\label{eq-T1Bound}
\ba{rcccl}
{1 \over \nu(1+t)} x(1) & \preceq_K & x(t) & \preceq_K & \nu \left(1 + {1 \over t} \right) x(1),\\
\\
{1 \over \nu(1+t)} s(1) & \preceq_{K^*} & s(t) & \preceq_{K^*} & \nu \left(1 + {1 \over t} \right) s(1).
\ea
\eeq
Substituting them to the bounds (\ref{eq-QBound}), (\ref{eq-GapUp}), and using the definition of $t_{\epsilon}$, we get
\beq\label{eq-QBound1}
\ba{rcl}
\la s(1), \bar x + x \ra + \la s, x(1) \ra & < &  {\nu(\nu + \epsilon)^2 \over \epsilon } \Big[\bar \tau + \la \bar s, x(1) \ra + \la s(1), \bar x \ra\Big],\\
\\
\bar \tau - \la c, x \ra + \la b, y \ra - \epsilon \tau & \leq  & (\nu + \epsilon)
\Big[ \bar \tau + \la s(1), \bar x \ra + \la \bar s, x(1) \ra\Big].
\ea
\eeq
Finally, by definition (\ref{def-Sigma}) of $\sigma \geq 1$, we have
\beq\label{eq-QBound2}
\ba{rcl}
\la \bar s, \bar x + x \ra + \la s, \bar x \ra & < &  {\nu \sigma (\nu + \epsilon)^2 \over \epsilon } \Big[\bar \tau + 2 \sigma \la \bar s,\bar x \ra\Big],\\
\\
\bar \tau - \la c, x \ra + \la b, y \ra - \epsilon \tau & \leq  & (\nu + \epsilon)
\Big[ \bar \tau + 2 \sigma \la \bar s, \bar x \ra \Big].
\ea
\eeq

Let us prove one more lemma, which depends on the following ``centrality measure'' for the reference points:
\beq\label{def-Mu}
\ba{rcl}
\mu & = & \max\limits_{\lambda} \Big\{ \lambda:  \; \nabla^2 F(\tilde x) \succeq \lambda \nabla^2 F(\bar x) \Big\}, \quad \bar \mu \Def \min\{ \mu, 1 \}.
\ea
\eeq
In view of definitions (\ref{eq-Rule}), it is easy to see that $\nabla^2 F_*(\tilde s) \succeq \mu \nabla^2 F_*(\bar s)$. For the choice of the reference points (\ref{eq-Choice}), we have $\tilde x = \bar x$ and $\mu = 1$.
\BL\label{lm-HessLow}
For any $u \in Q$, we have
\beq\label{eq-HessLow}
\ba{rcl}
\nabla^2 \Phi(u) & \succeq &  {\bar \mu \over c(K) \zeta^2(\epsilon)} \nabla^2 \Phi(0),
\ea
\eeq
where $\zeta(\epsilon) \Def  {\nu \sigma\over \epsilon} (\nu + \epsilon)^2  [1 + 2 \sigma \la \bar s, \bar x \ra ]$.
\EL
\proof
Let $u = (x,y,\tau) \in Q$. We extend it up to $\hat u = (x,s,y,\tau) $ with $s = \bar s + \tau c - A^*y$. 
In view of the choice (\ref{eq-Rule}), we have 
$$
\ba{rcl}
\la - \nabla F(\tilde x), \bar x + x \ra + \la s, - \nabla F^*(\tilde s) \ra & \refLEI{eq-QBound2}{1} & \zeta(\epsilon) \Def  {\nu \sigma\over \epsilon} (\nu + \epsilon)^2  [1 + 2 \sigma  \la \bar s, \bar x \ra ].
\ea
$$
In view of relation (\ref{eq-HLow1}) and definition (\ref{def-Mu}), we have
\beq\label{eq-FBar}
\ba{rcl}
\nabla^2 F(\bar x + x) & \succeq & {1 \over c(K) \zeta^2(\epsilon)} \nabla^2 F(\tilde x) \; \succeq \; {\bar \mu \over c(K) \zeta^2(\epsilon)} \nabla^2 F(\bar x) ,\\
\\
\nabla^2 F_*(s) & \succeq & {1 \over c(K) \zeta^2(\epsilon)} \nabla^2 F_*(\tilde s) \; \succeq \; {\bar \mu \over c(K) \zeta^2(\epsilon)} \nabla^2 F_*(\bar s).
\ea
\eeq

At the same time, 
\beq\label{eq-UGap}
\ba{rcl}
\ell(u) \Def \bar \tau - \la c, x \ra + \la b, y \ra - \epsilon \tau & \refLEI{eq-QBound2}{2}  & (\nu+\epsilon) [1 + 2 \sigma  \la \bar s, \bar x \ra] \; < \; \zeta(\epsilon). 
\ea
\eeq
Hence, we have ${\nabla \ell (\nabla \ell)^* \over \ell^2(u)} \succeq {1 \over \zeta^2(\epsilon)}\nabla \ell (\nabla \ell)^*$. Since $\ell(0) = 1$, this relation and (\ref{eq-FBar}) prove inequality  (\ref{eq-HessLow}).
\qed

We need to prove the last auxiliary statement.
\BL\label{lm-Zero}
Let $u \in Q$. Then, for any $v \in \E_u^*$ we have
\beq \label{eq-Zero}
\ba{rcl}
\lambda_u(v) & \leq & \sqrt{{c(K) \over \bar \mu} } \zeta(\epsilon) \lambda_0(v).
\ea
\eeq
\EL
\proof
Indeed, in view of representation (\ref{def-GNormR}) and relation (\ref{eq-HessLow}), we have
$$
\ba{rcl}
\lambda^2_u(v) & \leq & \max\limits_{h = (h_x,h_y,h_{\tau}) \in \E_u} \Big\{ 2 \la v, h \ra -  {\bar \mu \over c(K) \zeta^2(\epsilon)}\la \nabla^2 \Phi(0)h,h \ra: \; A h_x = h_{\tau} b \Big\}\\
\\
& = & {\bar \mu\over c(K) \zeta^2(\epsilon)} \lambda_0^2 \left( {c(K) \over \bar \mu} \zeta^2(\epsilon) v \right) \; = \; {c(K) \over \bar \mu}  \zeta^2(\epsilon) \lambda^2_0(v). \QR
\ea
$$

It remains to note that $\Phi(\cdot)$ is a $\bar \nu$-self-concordant barrier. Therefore,
$$
\ba{rcl}
\lambda^2_0(g) & \Def & \max\limits_{h = (h_x,h_y,h_{\tau}) \in \E_u} \Big\{ 2 \la g, h \ra - \la \nabla^2 \Phi(0)h,h \ra: \; A h_x = h_{\tau} b \Big\}\\
& = & \max\limits_{h = (h_x,h_y,h_{\tau}) \in \E_u} \Big\{ 2 \la \nabla \Phi(0), h \ra - \la \nabla^2 \Phi(0)h,h \ra: \; A h_x = h_{\tau} b \Big\}\\
& = & \max\limits_{h = (h_x,h_y,h_{\tau}) \in \E_u} \Big\{ 2 \la \nabla \Phi(0), h \ra - \la \nabla^2 \Phi(0)h,h \ra \Big\} \; \leq \; \bar \nu.
\ea
$$
Hence, $\lambda_{u_k}(g) \refLE{eq-Zero}  \sqrt{\bar \nu {c(K) \over \bar \mu} } \zeta(\epsilon) $.
Thus, the switching rule (\ref{eq-Switch}) provides us with the following sufficient condition for terminating the path-following scheme (\ref{met-SSPF}):
\beq\label{eq-Term}
\ba{rcl}
 \zeta(\epsilon) \sqrt{ \bar \nu \, {c(K) \over \bar \mu}}  & \leq & \tilde \beta t_k^{-1}.
\ea
\eeq
Since $t_k \leq \exp \left\{ - {\hat \gamma k \over \hat \beta + \sqrt{\nu}} \right\}$, we conclude that the number of iterations of the path-following scheme (\ref{met-SSPF}) cannot be bigger than
\beq\label{eq-CompPF}
\ba{c}
{\hat \beta + \sqrt{\nu} \over \hat \gamma} \ln \left( \zeta(\epsilon) \sqrt{ \bar \nu \, {c(K) \over \bar \mu}} \right).
\ea
\eeq

\section{Hot-Start Possibilities}\label{sc-Hot}
\SetEQ

Let us write down the optimality conditions defining the dual central path 
$$
\ba{rcl}
u_d(t)  & = &  (x_d(t), y_d(t), \tau_d(t)). 
\ea
$$
Denote by $v(t) \in \E_y$ the vector of optimal dual multipliers for the equality constraints in problem (\ref{def-HotCP}). Then we have $A x_d(t) = \tau_d(t) b$ and
\beq\label{eq-FOrder}
\ba{rcl}
t g_x + \nabla F(\bar x + x_d(t)) + {c \over \omega(t)} - A^* v(t) & = & 0,\\
\\
t g_y - A \nabla F^*(\bar s + \tau_d(t) c - A^* y_d(t)) - { b \over \omega(t)} & = & 0,\\
\\
t g_{\tau} + \la c, \nabla F^*(\bar s + \tau_d(t) c - A^* y_d(t)) + {\epsilon \over \omega(t)} + \la b, v(t) \ra & = & 0,
\ea
\eeq
where $\omega(t) = \bar \tau - \la c, x_d(t) \ra + \la b, y_d(t) \ra - \epsilon \tau_d(t)$. Using this trajectory, we can define the following {\em primal central path}:
\beq\label{def-PrimalCP}
\ba{rcl}
x_p(t) & = & - \omega(t) \nabla F^*(\bar s + \tau_d(t) c - A^* y_d(t)), \\
\\
s_p(t) & = & - \omega(t) \nabla F(\bar x + x_d(t)),\\
\\
y_p(t) & = & \omega(t) v(t).
\ea
\eeq
In view of (\ref{eq-FOrder}), it satisfies the following equations:
\beq\label{eq-Shadow}
\ba{rcl}
t  \omega(t) g_x  + c & = & s_p(t) + A^* y_p(t),\\
\\
t \omega(t) g_y  + A x_p(t) & = & b,\\
\\
t \omega(t) g_{\tau} + \epsilon   & = & \la c, x_p(t) \ra - \la b, y_p(t) \ra.
\ea
\eeq

Thus, the primal central path $u_p(t)$ inherits in a natural way the feasibility pattern of the point $\tilde w = (\tilde x, \tilde s, \tilde y)$.
\BL\label{lm-Feas}
\BI
\item
If $A \tilde x = b$, then $g_y = 0$ and $A x_p(t) = b$ for all $t > 0$.
\item
If $\tilde s + A^* \tilde y = c$, then $g_x = 0$ and $s_p(t) \Def c - A^*y_p(t) \in \inter K^*$ for all $t > 0$.
\EI
\EL

For example, if point $\tilde w$ belongs to the relative interior of the feasible set ${\cal F}$, then the primal central path connects it with the strictly feasible $\epsilon$-solution of our problem and this trajectory is always strictly feasible.

\section{Conclusion}\label{sc-Conc}
\SetEQ

In this paper, we proposed a new framework for finding a feasible $\epsilon$-solution of the primal-dual pair of Conic Problems (\ref{prob-PD}) by infeasible-start methods. Our methods are derived from one of the most important fact in Convex Analysis, an extended Gordan Theorem of Alternative.

For finding an approximate solution to the primal-dual pair of Comic Optimization problems,
these methods are applied to the problem (\ref{prob-Phi}), where the objective function is a self-concordant barrier subject to some linear equality constraints. It admits a unique interior optimal solution, allowing a significant acceleration in  the end of the process.

The problem (\ref{prob-Phi}) can be solved by all standard methods of Convex Optimization with linear equality constraints (Fast Gradient Methods, Quasi-Newton Methods, Second-Order Schemes, etc.). In this paper, we analyze only two polynomial-time short-step methods based on the theory of self-concordant functions: the Damped Newton Method and Path-Following scheme. Both methods end up with a local quadratic convergence in a small neighborhood of the solution.

Our path-following scheme follows the dual central path, which is defined as a trajectory of minimizers of a parametric barrier function. Using its gradient, we can define the primal central path. However, the situation is not symmetric since the primal path cannot be defined by an auxiliary potential.

In our approach, we can use the advantages of hot start. If the starting point is close to the solution, the path-following scheme needs less iterations to converge. Our method smoothly benefits from the feasibility level of the starting point and does not need to change from feasible to infeasible regime. Even for a feasible starting point, we have improved results since our trajectory goes directly to the optimal solution of the problem with no preliminary stages of the process.

Another advantage of our approach is the possibility to satisfy linear equality constraints with much higher accuracy than that one for the approximation of the optimal function value. This is achieved by the local quadratic convergence in the end of the process.

For the future research, it is very interesting to develop long-step variants of these methods and compare them with other unconstrained optimization schemes. It is a very rare situation when different approaches can directly compete on the same problem instance having very important practical applications.

\vspace{3ex}\noindent
{\bf Acknowledgment.} The author would like to thank Arkadi Nemirovski for very useful remarks and discussions of the results.

\section{Appendix: Complementary Facts on Conic\\  Optimization}\label{sc-App}
\SetEQ

Let $K \subset \E$ be a proper cone endowed with a $\nu$-regular barrier $F(\cdot)$. Let us fix a point $\bar x \in \inter K$ and define the set
\beq\label{def-Simplex}
\ba{rcl}
\Delta(\bar x) & = & \Big\{ x \in \inter K: \; \la \nabla F(\bar x), x - \bar x \ra \geq 0 \Big\}.
\ea
\eeq
\BL\label{lm-HLow}
For any $x \in \Delta(\bar x)$ we have
\beq\label{eq-HLow}
\ba{rcl}
\nabla^2 F(x) & \succeq & {1 \over c(K) \nu^2} \nabla^2 F(\bar x),
\ea
\eeq
where $c(K) = 4$. If the barrier $F(\cdot)$ is self-scaled, then we can take $c(K) = 2$.
\EL
\proof
Note that by the Theorem of Recession Direction, for any $x \in \Delta(\bar x)$, we have
\beq\label{eq-Size}
\ba{rcl}
\| x \|_{\bar x} & \refLE{eq-RecDir} & - \la \nabla F(\bar x), x \ra \; \refLE{def-Simplex} \; - \la \nabla F(\bar x), \bar x \ra \; \refEQI{eq-HomFX}{1} \; \nu.
\ea
\eeq

Consider now an arbitrary point $x \in \Delta(\bar x)$ and direction $u \in \E$. Note that 
$$
\ba{rcl}
W^F_1(x) & \Def & \{ y \in \E: \; \| y - x \|_{x} \leq 1\} \; \subseteq \; K. 
\ea
$$
Hence,  $x_u = x + {u \over \| u \|_{x}} \in K$. Without loss of generality, we can assume that 
$$
\ba{rcl}
\la \nabla F(\bar x), x_u - x \ra & = & \la \nabla F(\bar x), u \ra \;  \geq \; 0
\ea
$$
(otherwise, replace $u$ by $-u$). Thus, $x_u \in \Delta(\bar x)$. Note that for general cones, we have
$$
\ba{rcl}
\| x_u - x \|_{\bar x} \; \leq \; \| x \|_{\bar x} + \| x_u \|_{\bar x} & \refLE{eq-Size} & 2 \nu,
\ea
$$
and for self-scaled cones we have
$$
\ba{rcl}
\| x_u - x \|_{\bar x}^2 & = & \| x_u \|_{\bar x}^2 + \| x \|_{\bar x}^2 - 2 \la \nabla^2 F(\bar x) x_u, x \ra \; \leq \; \| x_u \|_{\bar x}^2 + \| x \|_{\bar x}^2 \; \refLE{eq-Size} \; 2 \nu^2.
\ea
$$
Thus, in both cases, we have $\| u \|_{\bar x}^2 \leq c(K) \nu^2 \| u \|^2_x$, where $c(K)$ is defined in the statement of Lemma \ref{lm-HLow}. And the later relation is exactly (\ref{eq-HLow}).
\qed

\BC\label{cor-HLow}
For all $x, \bar x \in \inter K$, we have
\beq\label{eq-HLow1}
\ba{rcl}
\nabla^2 F(x) & \succeq & {1 \over c(K) \la - \nabla F(\bar x), x \ra^2} \nabla^2 F(\bar x) .
\ea
\eeq
\EC
\proof
Let us choose $\tau = {\nu \over \la - \nabla F(\bar x), x \ra}$. Then $\la \nabla F(\bar x), \tau x \ra = - \nu \refEQI{eq-HomFX}{1} \la \nabla F(\bar x) , \bar x \ra$. Hence, $\tau x \in \Delta(\bar x)$ and therefore by (\ref{eq-HLow}) we have
$$
\ba{rcl}
{1 \over \tau^2} \nabla^2 F(x) \; \refEQI{eq-Hom12}{2} \; \nabla^2 F(\tau x) \; \succeq \; {1 \over c(K) \nu^2} \nabla^2 F(\bar x),
\ea
$$
and this is exactly (\ref{eq-HLow1}).
\qed

Let us justify now some properties of central path $(x(t), s(t), y(t))$ for problem (\ref{prob-PD}). For two positive values $t_1$, $t_2$, we have
$$
\ba{rcl}
\la s(t_1) - s(t_2), x(t_1) - x(t_2) \ra & = &0.
\ea
$$
Therefore,
\beq\label{eq-SXT1}
\ba{rcl}
\la s(t_1), x(t_2) \ra + \la s(t_2), x(t_1) \ra & = & \la s(t_1), x(t_1) \ra + \la s(t_2), x(t_2) \ra \\
\\
& \refEQI{eq-RCP}{1}& \nu \left[{1 \over t_1} + {1 \over t_2} \right].
\ea
\eeq
Since $W^{F_*}_1(s(t_1)) \subseteq K^*$, we have $s(t_1) - {s(t_2) \over \| s(t_2) \|_{s(t_1)}} \in K^*$. Therefore, in view of the Theorem of Recession Direction, we have
\beq\label{eq-SLow}
\ba{rcl}
s(t_1) & \succeq_{K^*} & {s(t_2) \over \| s(t_2) \|_{s(t_1)}} \; \stackrel{(\ref{eq-RecDir})}{\succeq_{K^*}}\; {s(t_2) \over \la - \nabla F_*(s(t_1)), s(t_2) \ra} \; \refEQI{eq-RCP}{1} \; {s(t_2) \over t_1 \la s(t_2), x(t_1) \ra} \; \stackrel{(\ref{eq-SXT1}) \atop \succeq_{K^*}} \; {t_2 s(t_2) \over \nu [t_1 + t_2]}.
\ea
\eeq
Similarly, for the primal central path, we get the following relation:
\beq\label{eq-XLow}
\ba{rcl}
x(t_1) & \succeq_{K} & {t_2 x(t_2) \over \nu [t_1 + t_2]}.
\ea
\eeq

\end{document}